\documentclass[12pt,oneside,reqno]{amsart}
\usepackage[all]{xy}
\usepackage{amsfonts,amsmath,oldgerm,amssymb,amscd}
\UseComputerModernTips
\numberwithin{equation}{section}
\usepackage[breaklinks]{hyperref}
\allowdisplaybreaks

\usepackage{enumerate}
\usepackage{caption} 
\newcommand{\Zstroke}{%
  \text{\ooalign{\hidewidth\raisebox{0.2ex}{--}\hidewidth\cr$Z$\cr}}%
}
\newtheorem{theorem}{Theorem}[section]

\newtheorem{lemma}[subsection]{{\bf Lemma}}
\newtheorem{coro}[subsection]{{\bf Corollary}}

\begin{document}

	\title[Solving Skolem problem for negative indexed $k$-generalized Pell numbers 
	]{Solving Skolem problem for negative indexed $k$-generalized Pell numbers
	} %\\ \today}

\author[M. Mohapatra]{M. Mohapatra}
\address{Monalisa Mohapatra, Department of Mathematics, National Institute of Technology Rourkela, Odisha-769 008, India}
\email{mmahapatra0212@gmail.com}

\author[P. K. Bhoi]{P. K.  Bhoi}
\address{Pritam Kumar Bhoi, Department of Mathematics, Institute of Mathematics $\&$Application, Andharua, Bhubaneswar, Odisha -751029, India}
\email{pritam.bhoi@gmail.com}

\author[G. K. Panda]{G. K. Panda}
\address{Gopal Krishna Panda, Department of Mathematics, National Institute of Technology Rourkela, Odisha-769 008, India}
\email{gkpanda\_nit@rediffmail.com}


\begin{abstract}
	In this
 paper, we address the Skolem problem for the $k$-generalized Pell sequence $(P_n^{(k)})_{n\geq2-k}$ extended to negative indices. We focus on identifying and bounding the indices $n<0$ for which $P_n^{(k)}=0.$ In particular, we establish that the zero multiplicity of $P_n^{(k)}$ is  $
\chi_k = \lfloor k^2/4\rfloor$
for all $k \in [4, 500].$
\end{abstract}

\maketitle
\pagenumbering{arabic}
\pagestyle{headings}
 \section{Introduction}
The $k$-generalized Pell sequence $(P_n^{(k)})_{n\geq2-k}$ extends the classical Pell sequence to a broader context for any integer $k \geq 2$. It is defined by the recurrence relation
$$P_n^{(k)}=2P_{n-1}^{(k)}+P_{n-2}^{(k)}+\cdots+P_{n-k}^{(k)}$$
for all $n\geq2$, with the initial terms $P_{n}^{(k)}=0$ for $n=-(k-2),-(k-3),\ldots,0,$ and $P_1^{(k)}=1$. Notably, when $k = 2$, the recurrence simplifies to that of the classical Pell sequence.

 One of the enduring unsolved questions in mathematics and computer science is whether a linear recurrence sequence ever reaches zero.
A classical formulation of this issue, known as the Skolem Problem, involves finding the set 
\begin{equation*}
\Zstroke(u) = \{n\in\mathbb{Z} : u_n = 0\}
\end{equation*}
for a given sequence $u=(u_n)_{n\in\mathbb{Z}}$.
If this set is finite, it's cardinality is referred to as the zero-multiplicity of the sequence $u$. However, there is no known general method or algorithm to compute the cardinality of $\Zstroke(u)$. A foundational result connected to this problem, established by Skolem, is the following:

\emph{{\bf Skolem Theorem}\cite{ShoreyTijdeman2008}, \cite{Skolem1934}{\bf:} If the coefficients of the linear recurrence sequence are rational, then the set $\Zstroke(u)$ is a union of finitely many arithmetic progressions together with a finite set.}

Subsequent extensions of this result were developed by Mahler, Lech, and others \cite{Lech1953, Mahler1935, MahlerCassels1956}, who demonstrated that when the coefficients are rational or algebraic, the set of zeros can be expressed as a finite combination of isolated points and arithmetic progressions. However, their proofs did not provide constructive methods. Later, Berstel and Mignotte \cite{BerstelMignotte1976} offered an effective approach specifically to identify these arithmetic progressions.

For non-degenerate linear recurrence sequences, the primary challenge is to precisely determine their zero multiplicity. Some advances have been achieved for sequences with integer coefficients and small order as discussed in \cite{Berstel1974}, \cite{Beukers1991}, \cite{BeukersTijdeman1984}, and \cite{Ward1959}.
For linear recurrence sequences of order $k\geq 4$, astronomical bounds depending on $k$ are known (e. g., Schmidt \cite{Schmidt1999} showed that $\Zstroke(u) \leq \exp(\exp(\exp(3k\log k)))$). For sequences with algebraic coefficients, van der Poorten and Schlickewei \cite{PoortenSchlickewei1991} derived bounds depending both on $k$ and the degree of the field $\mathbb{Q}(a_1,\ldots,a_k)$.
If the characteristic polynomial has real roots, Hagedorn \cite{Hagedorn2005} proved $\Zstroke(u) \leq 2k-3$. Garc\'{i}a et al. solved the Skolem problem for $k$-generalized Fibonacci numbers in \cite{ggl2020} and \cite{ggl2024}.

This paper focuses on analyzing the zero-multiplicity of the $k$-generalized Pell sequence $(P_n^{(k)})_{n \in \mathbb{Z}}$ extended to negative indices $n$ using the relation
\begin{equation*}
P_{n-k}^{(k)} = -P_{n-(k-1)}^{(k)} - P_{n-(k-2)}^{(k)} - \cdots -P_{n-2}^{(k)}- 2P_{n-1}^{(k)} + P_n^{(k)}.
\end{equation*}
Below are some specific values and ranges of negative indices $n \leq 0$ relevant to the zero-multiplicity of this sequence:
\begin{table}[h!]
\centering
\begin{tabular}{c|l|c}

$k$ & Indices  & Zero-Multiplicity \\ 
\hline 2 & $0$ & 1 \\
3 & $[-1,0]$ &2 \\
4 & $[-2, 0],-5$ & 4 \\
5 & $[-3, 0], [-7,-6]$ & 6 \\
6 &  $[-4, 0], [-9, -7]$, $-14$ & 9 \\
7 &  $[-5, 0], [-11, -8], [-17,-16]$ & 12\\
8 &  $[-6,0], [-13, -9], [-20,-18], -27$ & 16\\
9 &  $ [-7, 0], [-15, -10], [-23,-20],[-31,-30]$ & 20 \\
10 &  $ [-8, 0], [-17, -11], [-26,-22],[-35, -33],-44 $ & 25 \\
\end{tabular}

\label{tab:zero_multiplicity}
\end{table}
\begin{theorem}\label{main result 1}  
The largest nonnegative integer solution $n$ to the Diophantine equation $P_{-n}^{(k)} = 0$ satisfies the following bounds:
\begin{itemize}
    \item[(i)] If $k$ is even, then $n < 2k^{k^2}\log(16k^2).$ 
    \item[(ii)] If $k$ is odd, then $n < 7.5 \cdot 10^{14} \cdot 1.59^{k^3} \cdot k^{10} \cdot (\log k)^2.$ 
\end{itemize}
\end{theorem} 
We also derive the following corollary from this main result.
\begin{coro}\label{coro 1.2} Let $(\chi_k)_{k \geq 2}$ be the zero-multiplicity of the $k$-generalized Pell sequence. Then $\chi_2 = 1$, $\chi_3 = 2$ and 
$\chi_k =\lfloor k^2/4\rfloor,$ for $k\in[4,500].$ 
\end{coro}

\section{Auxiliary Results}A Baker-type lower bound for nonzero linear forms in logarithms of algebraic numbers is a useful tool for solving Diophantine equations. For such equations to be resolved efficiently, these lower bounds are essential. We will start by reviewing important definitions and significant results from the field of algebraic number theory.\\
Consider an algebraic number $\theta$ that has a minimal polynomial $$g(X) = c_0(X-\theta^{(1)})\cdot\cdot\cdot(X-\theta^{(k)}) \in \mathbb{Z}[X],$$where $c_0 > 0$ is the leading coefficient and $\theta^{(i)}$'s are conjugates of $\theta$. \\Accordingly, $$h(\theta) = \dfrac{1}{k}\biggr(\log c_0 + \sum\limits_{j=1}^{k}\max\{0,\log|\theta^{(j)}|\}\biggr),$$ is the $absolute$ $logarithmic$ $height$ of $\theta$. If $\theta = c/d$, $(d\neq 0)$ is a rational number with $c$ and $d$ being relatively prime, then $h(\theta)=\log$(max$\{|c|,d\}).$

We give some properties of the logarithmic height whose proofs can be found in \cite[Theorem B5]{b2018}. If $\theta$ and $\beta$ be two algebraic numbers, then 
\begin{equation*}
	\begin{split}
		& (i)\hspace{0.2cm} h(\theta\pm\beta) \leq h(\theta)+h(\beta)+\log 2,\\
		&(ii)\hspace{0.2cm} h(\theta\beta^{\pm 1})\leq h(\theta)+h(\beta),\\
		&(iii)\hspace{0.2cm} h(\theta ^k)=|k|h(\theta).
	\end{split}
\end{equation*}

Expanding on earlier notations, we introduce the following theorem from \cite{bgl2016} that enhances a result given by Matveev \cite{m2000}, as further developed by Bugeaud et al. \cite{bms2006}. This theorem provides a specific upper bound for our variables in the given equations.

\begin{theorem}\label{thm1}\cite{bgl2016}. Let $\theta_1,\ldots,\theta_t$ be positive real numbers in an algebraic number field $\mathbb{K}$ of degree $d_{\mathbb{K}}$ and $e_1, \ldots, e_t$ be nonzero integers. If $\Lambda = \prod\limits_{i=1}^{t} \theta_{i}^{e_i} -1$ is not zero, then 
	\begin{equation*}
	 \log |\Lambda| > -3 \cdot 30^{t+4} \cdot (t + 1)^{5.5} \cdot d_{\mathbb{K}}^2 (1 + \log d_{\mathbb{K}})(1 + \log tB)A_1 \cdots A_t,
	\end{equation*}
	where $B\geq \max\{|e_1|,\ldots,|e_t|\}$ and  $A_1, \cdots, A_t$ are positive integers such that $A_j \geq h'(\theta_j) = \max\{d_{\mathbb{K}}h(\theta_j),|\log \theta_j|, 0.16\},$ for $j= 1,\ldots,t$.
\end{theorem}
Next, we recall a refined version of the classical Baker--Davenport reduction lemma, due to Dujella and Peth\H{o} \cite{dp1998}. This lemma will be used to reduce the upper bounds on our variables.
\begin{lemma}\label{lem1}\cite{dp1998}. Let $M$ be a positive integer and $p/q$ denote a convergent of the continued fraction of the real number $\tau$ such that $q > 6M$. Consider the real numbers $A, B, \mu$ with $A > 0$ and $B > 1$. Let $\epsilon := \left\lVert \mu q\right\rVert - M\left\lVert \tau q\right\rVert$, where $\left\lVert .\right\rVert$ denotes the distance from the nearest integer. If $\epsilon > 0$, then there exists no solution to the inequality
	\begin{equation*}
		0 <|u\tau - v +\mu|<AB^{-w},
	\end{equation*}
	in positive integers $u,v,w$ with $u\leq M$ and $w\geq \dfrac{\log(Aq/\epsilon)}{\log B}$.
\end{lemma}
\begin{lemma}\label{lem2}\cite{sl2014}. Let $r\geq 1$ and $H > 0$ be such that $H > (4r^2)^r$ and $H>L/(\log L)^r$. Then $L < 2^rH(\log H)^r$.
\end{lemma} A linear recurrence sequence is called \emph{degenerate} if its characteristic polynomial has two distinct roots whose ratio is a root of unity, and \emph{non-degenerate} otherwise. Non-degenerate sequences are particularly significant in the study of zeros of linear recurrence sequences, since the absence of root-of-unity ratios ensures that the sequence does not vanish periodically. This property underlies the Skolem--Mahler--Lech Theorem, which describes the structure of the zero set of such sequences.  
A direct consequence of this theorem is the following result:
\begin{lemma} \cite{ggl2020}\label{lem dege}
For a non-degenerate linear recurrence sequence, the number of zero terms (i.e., its zero-multiplicity) is always finite.
\end{lemma} Before stating the next lemma, we recall the definition of the \emph{Mahler measure} of a polynomial. For
$
P(x) = a_d x^d + a_{d-1} x^{d-1} + \dots + a_0 \in \mathbb{Z}[x],
$ 
with roots \(\alpha_1, \dots, \alpha_d \in \mathbb{C}\), the Mahler measure is defined as
\[
M(P) = |a_d| \prod_{i=1}^{d} \max\{1, |\alpha_i|\}.
\]
\begin{lemma}\label{lem4}\cite{Dubickas2020}
Let $ P(x) \in \mathbb{Z}[x] $ be a polynomial of degree $ d \geq 2 $. Then, for any pair of its roots $ \eta, \zeta $ satisfying $ |\eta| \neq |\zeta| $, the following bounds hold:

\begin{itemize}
  \item[(i)] If both $ \eta $ and $ \zeta $ are nonreal, then
  \begin{equation}\label{nonreal}
  \big|\,|\eta| - |\zeta|\,\big| > \frac{\sqrt{3}}{2\left( \frac{d(d - 1)}{2} \right)^{\frac{d(d-1)}{4} + 1} M(P)^{\frac{d^3}{2}-d^2 - \frac{d}{2} + 1}}.
  \end{equation}
  
  \item[(ii)] If at least one of $ \eta $ or $ \zeta $ is real, then
  \begin{equation}\label{real}
  \big|\,|\eta| - |\zeta|\,\big| > \frac{1}{4d^{d^2/2 + d + 1} M(P)^{4d(d - 1) + 1}}.
  \end{equation}
  
  \item[(iii)] If both $ \eta $ and $ \zeta $ are real, then
  \begin{equation}\label{both real}
  \big|\,|\eta| - |\zeta|\,\big| \geq \frac{1}{2^{d^2/2 - 1} M(P)^{d - 1}}.
  \end{equation}
\end{itemize}
\end{lemma}
\section{\textit{k-}Generalized Pell Sequence}

In this section, we review several key properties of the 
$k$-generalized Pell sequence that will be useful in the subsequent analysis. The characteristic polynomial associated with this sequence is given by
$$
\Psi_{k}(x)=x^{k}-2 x^{k-1}-\cdots-x-1.
$$
In \cite{bhl2020}, Bravo et al. showed that 
\begin{equation}\label{bound gamma}
	\quad\varphi^{2}\left(1-\varphi^{-k}\right)<\gamma(k)<\varphi^{2}, \quad \text { for all } k \geq 2 , \quad \text{and} \quad \varphi=\dfrac{1+\sqrt{5}}{2}.
\end{equation}
Furthermore, in the same paper, they showed that $\Psi_{k}(x)$ is irreducible over $\mathbb{Q}[x]$ and has only one root $\gamma=\gamma(k)$ outside the unit circle. It is real and positive, so it satisfies $\gamma(k)>1$. The other roots are strictly inside the unit circle. Let $\gamma = \gamma_1, \gamma_2, \ldots, \gamma_k$ be the roots of the polynomial $\Psi_k(z)$, labeled such that
\begin{equation}\label{eq:rootsorder}
\gamma > |\gamma_2| \geq |\gamma_3| \geq \cdots \geq |\gamma_{k-1}| \geq |\gamma_k|, \quad \text{where } \gamma=\gamma(k).
\end{equation}
For $k \geq 2$, define
$$
	g_{k}(x):=\frac{x-1}{(k+1) x^{2}-3 k x+k-4}=\frac{x-1}{k\left(x^{2}-3 x+1\right)+x^{2}-1} . 
$$
In \cite{bh2020}, Bravo and Hererra proved that
\begin{equation}\label{eq2} \quad\left|g_{k}\left(\gamma_{i}\right)\right|<1,\quad2 \leq i \leq k.
\end{equation}
So, the number $g_{k}(\gamma)$ is not an algebraic integer. 
In addition, they showed that
\begin{equation}\label{binet like}
	P_{n}^{(k)}=\sum_{i=1}^{k} g_{k}\left(\gamma_{i}\right) \gamma_{i}^{{n-1}} \quad \text { and } \quad\left|P_{n}^{(k)}-g_{k}(\gamma) \gamma^{n-1}\right|<\frac{1}{2} 
\end{equation}
hold for all $n \geq 1$ and $k \geq 2$. In \cite{bhl2021}, it was shown that the logarithmic height of $g_k(\gamma)$ satisfies
\begin{equation}\label{eq8}  h(g_k(\gamma)) < 5 \log k \quad \text{for all } k \geq 2.  \end{equation}
A result established by Mignotte \cite{Mignotte1984} shows that for any pair of indices $1\leq i<j \leq k,$ the quotient $\gamma_i/\gamma_j$ is not a root of unity, so the $k$-generalized Pell sequence qualifies as a non-degenerate linear recurrence. Therefore, by Lemma \ref{lem dege}, it follows that the sequence has a finite number of zero terms, i.e., finite zero-multiplicity. 
The following lemma establishes some properties of $\gamma_i$ and $g_k(\gamma_i)$:
\begin{lemma}\label{lem:properties}
Let $\gamma = \gamma_1, \gamma_2, \ldots, \gamma_k$ be the roots of $\Psi_k(z)$ labeled such that their absolute values satisfy \eqref{eq:rootsorder}, and consider the function
\begin{equation*}
g_k(z) = \frac{z-1}{(k+1)(z^2+1)-3kz}.
\end{equation*}
Then:
\begin{itemize}
    \item[(i)] If $\gamma_i$ and $\gamma_j$ are two roots satisfying $|\gamma_i| > |\gamma_j|$, then
    \begin{equation*}
    \frac{|\gamma_i|}{|\gamma_j|} > 1 + 1.59^{-k^3}.
    \end{equation*}
    \item[(ii)] For all $k\geq 2$, $0.276 \leq g_k(\gamma) \leq 0.5.$
    \item[(iii)] For all $k\geq 2$ and $2\leq i\leq k$,
    $
    |g_k(\gamma_i)| < \min\left\{1, \dfrac{2}{k-2}\right\}.$
    \item[(iv)] For all $k\geq 2$,
    \begin{equation*}
    |\gamma_k| < 1 - \dfrac{\log\gamma}{2k}
    \quad \text{and} \quad
    |g_k(\gamma_k)| > \dfrac{\log\gamma}{2k(5k+4)}.
    \end{equation*}
    \item[(v)] If $ \gamma_i, \gamma_j $ are distinct roots of $ \Psi_k(z) $ such that $ |\gamma_i| = |\gamma_j| $, then $ \gamma_j = \overline{\gamma_i} $.
\end{itemize}
Note: Lemma 3.1(v) is a direct consequence of a paper by Mignotte \cite{Mignotte1984}, but again we
have given the proof of this result.
\end{lemma}

\textbf{Proof}: \begin{itemize}
 \item[(i)]    
{\bf Case 1: $ \gamma_i = \gamma$.}\\
In this case, we have
\begin{equation*}
\frac{|\gamma_i|}{|\gamma_j|} > |\gamma| > \frac{2d}{d + 1} \geq \frac{4}{3},
\end{equation*}
so the inequality holds trivially. (Of course, $ \gamma_j $ cannot be equal with $ \gamma $, since $ \gamma $ is the \vspace{0.3cm}conjugate of largest modulus.)\\ 
{\bf Case 2: At least one root is real.}\label{case 2}\\
Let $|\gamma_i|>|\gamma_j|$ and atleast one of $ \gamma_i $ or $ \gamma_j $ is real (if $ k $ is even). From Inequality \eqref{bound gamma}, we obtain an upper bound for $M(\Psi_k)$:
\[
M(\Psi_k)=\prod_{i=1}^{k}\max\{1,|\gamma_i|\}
= \gamma(k)<\varphi^{2}.
\] Substituting this upper bound into Inequality \eqref{real}, we get
\begin{equation*}\begin{split} 
|\gamma_i|-|\gamma_j|>\dfrac{1}{4k^{(k^2/2)+k+1}\cdot M(\Psi_k)^{4k(k-1)+1}}&>\dfrac{1}{4k^{(k^2/2)+k+1}\cdot(\varphi^2)^{4k(k-1)+1}}\\
&>\dfrac{1}{4k^{(k^2/2)+k+1}\cdot(2.7)^{4k(k-1)+1}}\\
&> \dfrac{1}{(10.8)k^{(k^2/2)+k+1}\cdot(54)^{k(k-1)}}.
\end{split}\end{equation*}

Since $(10.8)k^{(k^2/2)+k+1}\cdot(54)^{k(k-1)}<1.59^{k^3}$ for all $k\geq11$, it follows that $$\dfrac{1}{(10.8)k^{(k^2/2)+k+1}\cdot(54)^{k(k-1)}}>\dfrac{1}{1.59^{k^3}},$$ which implies $$\dfrac{|\gamma_i|}{|\gamma_j|}>1+1.59^{-k^3}.$$ For $k =
 2,4,6,8,10$ the negative real root and $k-2$ complex roots of the polynomial $\Psi_k$ can be readily computed using Python. In each of these cases, the required inequality can be directly verified through straightforward computation.\vspace{0.3cm}\\
{\bf Case 3: Both the roots are nonreal.} \\
Assume now that both $ \gamma_i $ and $ \gamma_j $ are nonreal algebraic numbers, with $ |\gamma_i| > |\gamma_j| $. Applying inequalities \eqref{nonreal} and \eqref{bound gamma}, we obtain
\begin{equation*}
|\gamma_i| - |\gamma_j| > \frac{\sqrt{3}}{(d(d - 1))^{d(d-1)/4 + 1} \cdot (\varphi^2)^{d^{3/2} - \frac{5d^2}{4} - \frac{d}{4} + 1}}.
\end{equation*}
Similar to the argument used in Case 2, and noting that $ |\gamma_j| < 1 $, it suffices to verify that
\begin{equation*}
(\varphi^2)^{d^{3/2} - \frac{5d^2}{4} - \frac{d}{4} + 1} (d(d - 1))^{(d(d-1)/4) + 1} < \sqrt{3} \cdot 1.59^{d^3}.
\end{equation*}
This inequality can be confirmed using Python for all $ d \geq 2 ,$ completing the argument.
 \item[(ii)] This is proved in \cite{bh2020}.
 \item[(iii)] For $k\geq 4$, the statement follows from \cite{bh2020}, and for $k = 2,3,$ it can be easily verified through direct computation that $|g_k(\gamma_i)| < 1.$
\item[(iv)] This proof is followed from \cite{ggl2020} which uses Mahler measure bounds (see \cite{Waldschmidt2000}) and estimates on the roots. Denoting the Mahler measure of $\gamma$ by $M(\gamma)$ and the house (maximum absolute value among conjugates) by $\overline{|\gamma|},$ we obtain from \eqref{eq:rootsorder}:
\begin{equation*}
|1/\gamma_k|^k = \overline{|1/\gamma_k|}^k \geq M(1/\gamma_k) = M(\gamma_k) = \gamma.
\end{equation*}
Thus,
\begin{equation*}
|\gamma_k| \leq \gamma^{-1/k} = \exp\left( -\frac{\log\gamma}{k} \right) < 1 - \frac{\log\gamma}{2k},
\end{equation*}
where the last inequality uses $\exp(-x) < 1 - \frac{x}{2}$ for small positive $x$.
Now, since
\begin{equation*}
\left|k(\gamma_i^2-3\gamma_i+1)+\gamma_i^2-1\right| =\left|(k+1)(\gamma_i^2+1)-3k\gamma_i-2\right| < 5k+4,
\end{equation*}
and
$
|\gamma_k-1| > (\log\gamma)/2k,
$
we conclude
\begin{equation*}
|g_k(\gamma_k)| > \frac{\log\gamma}{2k(5k+4)}.
\end{equation*}
\item[(v)] We have done this proof by following the approach in \cite{p2020}. Define 
$$
\delta_k(x) = (x - 1)\Psi_{k}(x) = x^{k+1} - 3x^k +x^{k-1}+ 1.
$$
Apart from 1, the polynomial $ \delta_k(x) $ shares the same roots as $ \Psi_{k}(x) $. Moreover, all roots of $ \Psi_{k}(x) $ are simple (i.e., have multiplicity one).
If both $ \gamma_i, \gamma_j $ are real numbers and satisfy $ |\gamma_i| = |\gamma_j| ,$ then it must be that $ \gamma_j = -\gamma_i $. Assuming $ \gamma_i > 0 $, we consider the expression
$$
\gamma_i^{k+1} - 3\gamma_i^k+\gamma_i^{k-1} + 1 = (-1)^{k+1}\gamma_i^{k+1} - 3(-1)^k\gamma_i^k+(-1)^{k-1}\gamma_i^{k-1} + 1 = 0.
$$
Simplifying this leads to
$
\gamma_i^2-3\gamma_i+1  = (-1)^{k+1}(\gamma_i^2+3\gamma_i+1).
$
If $ k $ is odd, then $-\gamma_i=\gamma_i,$ implies $\gamma_i=0,$ a contradiction. For $ k $ is even, the equation would imply $ \gamma_i^2 = -1 $, which is also impossible.
Now, suppose that $ \gamma_i$ and $\gamma_j $ are non-real complex numbers. Let their common absolute value be $ r $ and their arguments be $ \omega_1$ and $\omega_2 $, respectively, so that
$$
\gamma_i = r(\cos\omega_1 + i\sin\omega_1), \quad \gamma_j = r(\cos\omega_2 + i\sin\omega_2).
$$
Since $ \delta_k(\gamma_i) = \delta_k(\gamma_j) = 0 ,$ we have
$$
r^{k+1} l_1 - 3r^k m_1 +r^{k-1}n_1+ 1 = 0 = r^{k+1} l_2 - 3r^k m_2 +r^{k-1}n_2+ 1 ,
$$
with
$$
\begin{aligned}
l_1 &= \cos((k+1)\omega_1) + i\sin((k+1)\omega_1), \\
l_2 &= \cos((k+1)\omega_2) + i\sin((k+1)\omega_2), \\
m_1 &= \cos(k\omega_1) + i\sin(k\omega_1), \\
m_2 &= \cos(k\omega_2) + i\sin(k\omega_2),\\
n_1 &= \cos((k-1)\omega_1) + i\sin((k-1)\omega_1), \\
n_2 &= \cos((k-1)\omega_2) + i\sin((k-1)\omega_2). 
\end{aligned}
$$
After some simplification, we obtain
\begin{equation}\label{eqn}
r^2(l_1 - l_2)+(n_1-n_2) = 3r(m_1 - m_2).
\end{equation}
Here \begin{equation*}
    \begin{split}
\arg(l_1 - l_2) &= \frac{\pi}{2} + \frac{(k+1)(\omega_1 + \omega_2)}{2}, \\
\arg(m_1 - m_2) &= \frac{\pi}{2} + \frac{k(\omega_1 + \omega_2)}{2},\quad\text{and }\\
\arg(n_1 - n_2) &= \frac{\pi}{2} + \frac{(k-1)(\omega_1 + \omega_2)}{2}.
   \end{split}
\end{equation*}
Let's denote their corresponding moduli by
\[
s = |l_1 - l_2|, \quad t = |m_1 - m_2|, \quad\text{and } \quad u = |n_1 - n_2|.
\]
We have $ s = 0 $ if and only if $ l_1 = l_2 $. Then it leads to $(n_1-n_2) = 3r(m_1 - m_2).$ Thus, $n_{1}-n_{2}$ is a real scalar multiple of $m_{1}-m_{2}$ 
and the difference of their arguments is
\[
\arg(m_1 - m_2) - \arg(n_1 - n_2) = \frac{\omega_1 + \omega_2}{2}.
\]
This shows that the vectors $ m_1-m_2$ and $ n_1-n_2$ point in different directions in the complex plane and cannot satisfy a real scalar multiple relation unless $\omega_1 + \omega_2 \equiv 0 \pmod{2\pi}$. However, $\omega_1 + \omega_2 \equiv 0 \pmod{2\pi}$ implies $l_1 \neq l_2$,
contradicting our assumption. Consequently, the only way the equation can hold is if
\[
n_1 - n_2 = 0 \quad \text{and} \quad m_1 - m_2 = 0,
\]
which together with $l_1 = l_2$ imply $ e^{i\omega_1} = e^{i\omega_2} $ by Euler's identity, hence $ \omega_1 = \omega_2 $, i.e., $ \gamma_i = \gamma_j $, a contradiction. In the sequel, we may assume  $ s,t,u \neq 0 $. Now Equation \eqref{eqn} implies
\begin{equation*}
 r^2s \cdot e^{-i(k+1)\frac{\omega_1 + \omega_2}{2}}+u\cdot e^{-i(k-1)\frac{\omega_1 + \omega_2}{2}} = 3rt \cdot e^{-ik\frac{\omega_1 + \omega_2}{2}}, \end{equation*}which on simplification yields 
\begin{equation}\label{3.7}
    r^2s \cdot e^{-i\frac{\omega_1 + \omega_2}{2}}+u \cdot e^{i\frac{\omega_1 + \omega_2}{2}}= 3rt.
\end{equation}
Comparing both sides of Equation \eqref{3.7}, we get $$(r^2s+u)\cos\biggr({\dfrac{\omega_1+\omega_2}{2}}\biggr)=3rt \quad\text{and} \quad(u-r^2s)\sin\biggr({\dfrac{\omega_1+\omega_2}{2}}\biggr)=0.$$ Thus, either \(\sin((\omega_1+\omega_2)/2)=0\) (which immediately gives \(\omega_1+\omega_2\equiv 0\pmod{2\pi}\)), or $r^2 s = u$. Using the polar form of complex numbers, we obtain
\[
s = |l_1 - l_2| = (l_1 - l_2){e^{-i \arg(l_1-l_2)}}, \quad
u = |n_1 - n_2| = (n_1 - n_2){e^{-i \arg(n_1-n_2)}}.
\]
Then the equality $r^2 s = u$ becomes
\[
r^2 (l_1 - l_2) e^{i (\arg(n_1-n_2) - \arg(l_1-l_2))} = n_1 - n_2.
\]
Since $r^2>0$, the left-hand side is a positive scalar multiple of $l_1-l_2$ rotated by $\arg(n_1-n_2) - \arg(l_1-l_2)$. For equality of nonzero complex numbers, the directions must coincide:
$$
\arg(n_1-n_2) - \arg(l_1-l_2) \equiv 0 \pmod{2\pi}.
$$
But
$
\arg(n_1-n_2) - \arg(l_1-l_2) = -(\omega_1+\omega_2)/2.
$
So, the equality can only hold if $\omega_1+\omega_2 \equiv 0 \pmod{2\pi}.$  
Therefore, in both the cases we get $\omega_2 = -\omega_1$, i.e., $\gamma_j = \overline{\gamma_i}$, which is the desired conclusion.  
\end{itemize}

\begin{lemma}\label{lem2.2 ref}
Let $\gamma_1,\dots,\gamma_{2k}$ be the roots of $\Psi_{2k}(x)$. Then
$\gamma_{2k}$ is real and negative with
$
|\gamma_{2k}|<|\gamma_j|$ for $j=1,\ldots,2k-1.$
\end{lemma}

\textbf{Proof}:
Let
$
T_{2k}(x) = x^{2k+1} \delta_{2k}(-1/x) 
= x^{2k+1} - x^2 - 3x - 1.
$
Since 
\[
T_{2k}(1) = 1 - 1 - 3 - 1 = -4 < 0, \qquad 
T_{2k}(2) = 2^{2k+1} - 11 > 0,
\]
it follows that $T_{2k}$ has a real root $\lambda \in (1,2)$. 
As $T_{2k}$ is strictly increasing on $(1,\infty)$, $\lambda$ is its only real root in this interval. Now fix $\mu > \lambda$. Then $T_{2k}(\mu) > 0$. 
If $z \in \mathbb{C}$ with $|z|=\mu$, we have
\[
|z^2+3z+1| \le \mu^2 + 3\mu + 1 < \mu^{2k+1}.
\]
Hence, on the circle $|z|=\mu$, one has 
$|z^{2k+1}| > |-(z^2+3z+1)|$. 
By Rouché’s theorem \cite{f1990}, $T_{2k}(z)$ and $z^{2k+1}$ have the same number of zeros in $|z|<\mu$, namely $2k+1$. 
As $\mu>\lambda$ is arbitrary, all roots of $T_{2k}$ lie in the closed
disc $|z|\le \lambda$, with equality holds only for $z=\lambda$. 
In view of (v) of Lemma \ref{lem:properties}, $\lambda$ is the only root of $T_{2k}$ with absolute value $\lambda$, thus the absolute value of all other roots must be smaller, i.e., $\lambda=-1/\gamma_{2k}.$

\begin{lemma}\label{lem2.2}
Let $k$ be even. Then
\begin{equation*} \frac{|\gamma_{k-1}|}{|\gamma_k|} > 1 + \frac{1}{k^{k^2}}. \end{equation*}
\end{lemma}
\textbf{Proof}: Since $\gamma_k$ is real (by Lemma \ref{lem2.2 ref}), it follows from Equation \eqref{real} that
\begin{equation*} ||\gamma_{k-1}| - |\gamma_k|| > \frac{1}{4k^{k^{\frac{k^2}{2} + k + 1}}(\varphi^2)^{(k(k-1))+1}}. \end{equation*}
Thus, we have
\begin{equation*} \frac{|\gamma_{k-1}|}{|\gamma_k|} > 1 + \frac{1}{4|\gamma_k|k^{\frac{k^2}{2} + k + 1}(\varphi^2)^{(4{k(k-1))+1}}}. \end{equation*}
Using (iv) of Lemma \ref{lem:properties} and the fact that $\gamma>\varphi$, we have $$|\gamma_k| < 1 - \frac{\log \gamma}{2k} < 1 - \frac{\log\varphi}{2k}.$$ Thus, it suffices that
\begin{equation*} \frac{1}{4(1 - (\log(\varphi)/2k))k^{\frac{k^2}{2} + k + 1}(\varphi^2)^{(4{k(k-1))+1}}} > \frac{1}{k^{k^2}}, \end{equation*}
which holds true for every $k \geq 2$. $\hfill \Box$

\begin{lemma}\label{lem 3.4}
If $k$ is odd, then
$
\Psi_k(z)
$ has exactly one positive real root, and all remaining roots occur in complex conjugate pairs.
\end{lemma}
\textbf{Proof}:
As observed in \cite{bhl2020}, $\Psi_k(z)$ has exactly one real positive root.
Next, we show that there is no real negative root of $\Psi_k(z)$. Substituting $x=-y$ with $y>0$ gives
\[
\Psi_k(-y) = (-y)^k - 2(-y)^{k-1} - \sum_{i=0}^{k-2} (-y)^i.
\]
Since $k$ is odd,
\[
\sum_{i=0}^{k-2} (-y)^i = \frac{1 - (-y)^{k-1}}{1+y} = \frac{1 - y^{k-1}}{1+y}.
\]
Therefore,
\[
\Psi_k(-y) = -y^k - 2y^{k-1} - \frac{1 - y^{k-1}}{1+y} - 1.
\]
Multiplying both sides by $1+y>0$ yields
\[
(1+y)\Psi_k(-y) = -y^{k+1} - 3y^k - y^{k-1} - y - 2.
\]
Every term on the right-hand side is negative for $y>0$, hence $(1+y)\Psi_k(-y)<0$. Thus, $\Psi_k(-y)<0$ for all $y>0$, showing that $\Psi_k(z)$ has no negative real root. Combining this with the existence of a unique positive real root and the fact that
non-real roots occur in conjugate pairs yields that the remaining $k-1$ roots
are complex conjugate pairs.

\section{Proof of Theorem \ref{main result 1}}

Assume that $k\geq 4$. We consider the Diophantine equation $P_{-n}^{(k)} = 0$.  
By Equation \eqref{binet like}, this problem reduces to solving the exponential Diophantine equation 
\begin{equation}\label{eqn1}
\sum_{i=1}^{k} g_k(\gamma_i)\gamma_i^{-(n+1)} = 0
\end{equation}
for integers $n\geq 0$,
where $\gamma_1, \ldots, \gamma_k$ are the roots of the characteristic polynomial $\Psi_k(z)$, ordered as in \eqref{eq:rootsorder}.
\vspace{-0.7cm} We proceed by treating the cases of even and odd $k$ separately.\\
\subsection*{Case 1: \textit{k} is even}
Since $\Psi_k(z)$ possesses two real roots with opposite signs when $k$ is even, it follows that $\gamma_k$ is a negative real root.  
Equation \eqref{eqn1} can be rewritten as:
\begin{equation*}
g_k(\gamma_k)\gamma_k^{-(n+1)} = -\sum_{j=1}^{k-1} g_k(\gamma_j)\gamma_j^{-(n+1)}.
\end{equation*}
Taking absolute values, we get
\begin{equation*}
\left| g_k(\gamma_k)\gamma_k^{-(n+1)} \right| \leq \sum_{j=1}^{k-1} \left| g_k(\gamma_j)\gamma_j^{-(n+1)} \right|.
\end{equation*}
Estimates (iii) and (iv) of Lemma \ref{lem:properties} yield
\begin{equation*}
\begin{split}
\left( \frac{\log \gamma}{2k(5k+4)} \right) |\gamma_k|^{-(n+1)} &\leq |g_k(\gamma_k)|\gamma_{k}|^{-(n+1)}\\&\leq\biggr(\sum_{j=1}^{k-1}|g(\gamma_j)|\biggr)|\gamma_{k-1}|^{-(n+1)}\\&<0.5||\gamma_{k-1}|^{-(n+1)}|.
\end{split}
\end{equation*}
Applying Lemma \ref{lem2.2}, we derive
\begin{equation}\label{11}
\left(1+ \frac{1}{k^{k^2}} \right)^{n+1}<\left( \frac{|\gamma_{k-1}|}{|\gamma_k|} \right)^{n+1} < \dfrac{4k(5k+4)}{\log \gamma} < 41k^2.
\end{equation}
Now, taking logarithms on both sides of the above inequality, we find
\begin{equation*}
\dfrac{n+1}{2k^{k^2}}<(n+1)\log\left(1+\frac{1}{k^{k^2}}\right) < \log(41k^2),
\end{equation*}
which implies $
n < 2k^{k^2}\log(41k^2).$
\subsection*{Case 2: \textit{k} odd}
 In this case, $\Psi_k(z)$ has only one real positive root $\gamma$, and the rest are complex conjugates (by Lemma \ref{lem 3.4}). For linear recurrence sequences, Mignotte~\cite{Mignotte1975} showed that the relation $P^k_{-n}=0$ implies $n<c(k)$ for some effectively computable
constant depending only on $k$. However, the
value of $c(k)$ was not made explicit. Here, we determine $c(k)$ explicitly.
Define
\begin{equation*}
\Lambda=1 + \frac{g_k(\gamma_{k-1})}{g_k(\gamma_k)}\left( \frac{\gamma_k}{\gamma_{k-1}} \right)^{n+1} \quad\text{  for all }k\geq2 .
\end{equation*}
In view of Equation \eqref{eq:rootsorder} have that $\gamma^{-1}<|\gamma_2^{-1}
  |\leq\cdots\leq|\gamma_k^{-1}|$. Furthermore, since $k$ is odd, we have
 $|\gamma_{k}|^{-1}=|\gamma_{k-1}|^{-1}>|\gamma_{k-2}|^{-1}.$ Thus, we separate the terms on the left-hand side of Equation \eqref{eqn1}
 that involve $\gamma_k$ and $\gamma_{k-1}$. To proceed further, we divide both sides of the equation by $g_k(\gamma_k)\gamma_k^{-(n+1)}$
 and then take the absolute value. This yields
\begin{equation*}
|\Lambda|=\biggr|\sum_{i=1}^{k-2} \frac{g_k(\gamma_i)}{g_k(\gamma_k)}\left( \frac{\gamma_k}{\gamma_i} \right)^{n+1}\biggr|.
\end{equation*}
Applying estimates (ii),
 (iii) from Lemma \ref{lem:properties}, we see that
\begin{equation}\label{eqn2}
\sum_{i=1}^{k-2} \left| \frac{g_k(\gamma_i)}{g_k(\gamma_k)} \right| \left| \frac{\gamma_k}{\gamma_i} \right|^{n+1} < \frac{0.5(k-2)}{(k-2)g_k(\gamma_k)} \left|\frac{\gamma_k}{\gamma_{k-2}}\right|^{n+1}.
\end{equation}
Since $|\gamma_{k-2}|>|\gamma_k|,$ it follows from (i) and (iv) of Lemma \ref{lem:properties} that
\begin{equation}\label{eqn3}
  |\Lambda|  <\dfrac{2k(5k+4)}{(1+1.59^{-k^3})^{n+1}}.
\end{equation}
Apply Matveev's theorem (Theorem \ref{thm1}) with:
\begin{align*}
    &\theta_1 = -\frac{g_k(\gamma_{k-1})}{g_k(\gamma_k)}, \quad \theta_2 = \frac{\gamma_k}{\gamma_{k-1}}, \\
    &e_1 = 1, \quad e_2 = n + 1.
\end{align*}
We consider the number field $\mathbb{K} = \mathbb{Q}(\gamma_k, \gamma_{k-1})$ which has degree at most $[\mathbb{K} : \mathbb{Q}]=d_{\mathbb{K}} \leq k^2$. Using the
 properties \eqref{eq:rootsorder} and \eqref{eq8}, we get
 $$h(\theta_1)<10\log k \quad\text{ and } \quad h(\theta_2)<(\log 6)/k.$$
Moreover, since $g_k(\gamma_k)$ and $g_k(\gamma_{k-1})$ are algebraic conjugates (and also complex conjugates), their quotient $\theta_i$ lies on the unit circle, i.e., $|\theta_i| = 1$ and $|\log \theta_i| < \pi$, for $i = 1, 2$. Hence, we can take
\begin{equation*}A_1 = 10k^2 \log k, \quad A_2 = 1.8k, \quad B = n + 1.
\end{equation*}
We proceed under the assumption that $\Lambda \neq 0$. To the contrary, suppose $\Lambda= 0$. Then it would follow that
\begin{equation}\label{eqn4}\left(\frac{\gamma_k}{\gamma_{k-1}}\right)^{n+1} = -\left(\frac{g_k(\gamma_k)}{g_k(\gamma_{k-1})}\right). \end{equation}
 Now, consider a field automorphism $\sigma$ in the Galois group of the splitting field of $\Psi_k(z)$ over $\mathbb{Q}$ such that $\sigma(\gamma_k) = \gamma.$ Let $\gamma_j = \sigma(\gamma_{k-1})$. Clearly, $j \neq 1$. Then $\sigma(g_k(\gamma_k)) = g_k(\gamma)$ and $\sigma(g_k(\gamma_{k-1})) = g_k(\gamma_j)$. Conjugating relation \eqref{eqn4} by $\sigma$ and taking absolute values, we get
\begin{equation}\label{4.6}\left|\frac{\gamma}{\gamma_j}\right|^{n+1} = \left|\frac{g_k(\gamma)}{g_k(\gamma_j)}\right|.
\end{equation}
 Since $(g_k(\gamma_{}))<(g_k(\gamma_{j}))$, where $\gamma_j<\gamma_{}$ (By Lemma 3.2 of \cite{bhl2020}), the L.H.S of Equation \eqref{4.6} is $>1$ but the R.H.S is $<1$, which is a contradiction. Thus, we may assume that $\Lambda \neq 0$.
By Theorem \ref{thm1}, we get
\begin{equation}\label{eqn5}\begin{split}
\log |\Lambda| &> -1.7 \cdot 10^{13} \cdot k^7(1 + \log k^2)(1 + \log(2n + 2)) \log k\\ &> -2.1 \cdot 10^{14} \cdot k^7 \log(n + 1)(\log k)^2, \end{split}\end{equation}
where we have used the fact that $(1+\log k^2)<3.5 \log k, \text{ for all }k\geq2$ and $(1+\log(2n+2))<3.5 \log(n+1).$
Combining Equations \eqref{eqn3} and \eqref{eqn5}, we obtain
\begin{equation*}-2.1 \cdot 10^{14} \cdot k^7 \log(n + 1)(\log k)^2 < \log(2k(5k + 2)) - (n + 1) \log(1 + 1.59^{-k^3}).
\end{equation*}
Since the inequalities
\begin{equation*}\log(2k(5k + 2)) < 3.5 \log k \quad \text{and} \quad \log(1 + 1.59^{-k^3}) > \left(1 + 1.59^{k^3}\right)^{-1}
\end{equation*}
hold for all $k \geq 5$, we have
\begin{equation}\label{eqn6}\frac{n + 1}{\log(n + 1)} < 3.1 \cdot 10^{14} \cdot 1.59^{k^3} \cdot k^7 \cdot (\log k)^2. \end{equation}
Using Lemma \ref{lem2} in Inequality \eqref{eqn6}, we obtain an upper bound of the form
\begin{equation*}
\begin{split}n &<c(k)= 2 \cdot \left(3.1\cdot 10^{14} \cdot 1.59^{k^3} \cdot k^7 \cdot (\log k)^2\right) \cdot \log\left(3.1 \cdot 10^{14} \cdot 1.59^{k^3} \cdot k^7 \cdot (\log k)^2\right)\\ &< 7.5 \cdot 10^{14} \cdot 1.59^{k^3} \cdot k^{10} \cdot (\log k)^2,
\end{split}
\end{equation*}
where we have used
\begin{equation*}\log\left(3.1\cdot 10^{14} \cdot 1.59^{k^3} \cdot k^7 \cdot (\log k)^2\right) < 1.2k^3 \quad \text{for all } k \geq 4.
\end{equation*}
This completes the proof of part (ii) of Theorem \ref{main result 1}.\\
\section{Proof of Corollary \ref{coro 1.2}}

In this section, we determine the zero-multiplicity of \( (P_n^{(k)})_{n \in \mathbb{Z}} \), for each $k$ in the range \( [4, 500] \). Rather than directly applying the explicit upper bounds given in Theorem \ref{main result 1}, we adopt the approach used in its proof to derive more refined and practical bounds on $n $ for these values of $k$.
Our analysis considers separate cases based on whether \( k \) is even or odd.
\subsection*{Case 1: \textit{k} is even}
By analyzing the inequality derived in Theorem \ref{main result 1}, specifically equation \eqref{11}, we arrive at
\[
\left( \frac{|\gamma_{k-1}|}{|\gamma_k|} \right)^{n+1} < 41k^2.
\]
Applying logarithms to both sides, we get $$ (n+1)\log\left( \frac{|\gamma_{k-1}|}{|\gamma_k|} \right)<\log(41k^2).$$
For each even integer $k$ in the range \( [4, 500] \), we used this inequality to compute an explicit upper bound on \( n \). The outcome of this computation is formalized in the following lemma:

\begin{lemma}\label{lem4.1}
Suppose \( k \) be even, and \( L_k \) represent an upper bound for those integers \( n_k \geq 0 \) satisfying \( P^{(k)}_{-n_k} = 0 \). Then, for all \( k \in [4, 500] \),
\[
L_k \in [154, 9131366].
\]
\end{lemma}
\subsection*{Case 2: \textit{k} is odd}
Using Equations \eqref{eqn2} and \eqref{eqn5}, we obtain:
\[
2.1 \cdot 10^{14} \cdot k^7 \log(n + 1)(\log k)^2 > \log\left|\frac{2g_k(\gamma_k)}{3}\right| + (n + 1) \log\left|\frac{\gamma_{k-2}}{\gamma_k}\right|.
\]
Through a more refined analysis, we arrive at the following result:
\begin{equation}\label{eqn8}
\text{If } P^{(k)}_{-n} = 0 \text{ and } k \in [4, 500], \text{ then } 0 \leq n < 2 \cdot 10^{48}. 
\end{equation}
Using the identity for complex logarithms, namely \( \log w = \log|w| + i \arg w \), we recall the power series expansion
\[
\log (1+w) = \sum_{n=1}^\infty(-1)^{n-1}\dfrac{w^n}{n}\quad \text{valid for all } w\in \mathbb{C} \quad\text{ with } |w|<1.
\]
This series implies the estimate $|\log(1+w)|\leq 2|w|$ whenever $|w|\leq1/2.$ Returning to Equation \eqref{eqn2}, we carried out numerical checks confirming that the condition $|\Lambda| < 1/2$ holds for all odd $k
 \in [4,500]$ under the assumption $n > k^3.$ This assumption is justified as the bounds ultimately derived for $n$ begin well above $k^3$ in each instance. It is important to remember that the complex logarithm behaves additively only modulo $2\pi i$. Taking this into account, Equation \eqref{eqn2} leads to the conclusion:
 \begin{equation} \label{eqn9}
 \begin{split}
0 < |\Lambda| &:=
\left| \log\left( -\frac{g_k(\gamma_{k-1})}{g_k(\gamma_k)} \right) + (n+1) \log\left( \frac{\gamma_k}{\gamma_{k-1}} \right) + 2\pi i \ell \right|\\
&\leq \dfrac{0.5} {\left|\gamma_{k-2}/\gamma_k\right|^{n+1}} \cdot \left| \frac{1}{g_k(\gamma_k)} \right|. 
\end{split}
 \end{equation}
Using \( \gamma_k = \overline{\gamma_{k-1}} \), we have
\[
\log\left( \frac{\gamma_k}{\gamma_{k-1}} \right) = \log\left( \frac{\gamma_k}{\overline{\gamma_k}} \right) = \log\left( e^{2i \arg(\gamma_k)} \right) = 2i \arg(\gamma_k).
\]
Also,
\[
\log\left( -\frac{g_k(\gamma_{k-1})}{g_k(\gamma_k)} \right)
= \log\left( -\frac{\overline{g_k(\gamma_k)}}{g_k(\gamma_k)} \right)
= \log(-1) + \log\left( \frac{\overline{g_k(\gamma_k)}}{g_k(\gamma_k)} \right)
= i\pi - 2i \arg(g_k(\gamma_k)).
\]
Because \( \gamma_k = \overline{\gamma_{k-1}} \), we can rewrite
\[
\log\left( \frac{\gamma_k}{\gamma_{k-1}} \right) = 2i \arg(\gamma_k), \quad 
\log\left( -\frac{g_k(\gamma_{k-1})}{g_k(\gamma_k)} \right) = i(\pi - 2\arg(g_k(\gamma_k))).
\]
Thus, Inequality \eqref{eqn9} becomes
\begin{equation}\label{eqn10}
0 <|\Lambda|/\pi= \left| \left( -\frac{2 \arg(\gamma_k)}{\pi} \right)(n+1) - (2\ell + 1) + \frac{2\arg(g_k(\gamma_k))}{\pi} \right|
< 0.2 \left| g_k(\gamma_k) \right|^{-1} \left| \frac{\gamma_{k-2}}{\gamma_k} \right|^{-(n+1)}. 
\end{equation}
Now define
\[
\tau_k := -\frac{2 \arg(\gamma_k)}{\pi}, \quad
\mu_k := \frac{2 \arg(g_k(\gamma_k))}{\pi}, \quad
A_k :=  |g_k(\gamma_k)|^{-1}, \quad
B_k := \left| \frac{\gamma_{k-2}}{\gamma_k} \right|.
\]
Then Inequality \eqref{eqn10} takes the form
\begin{equation}\label{eqn21}
0 < |u \tau_k - v + \mu_k| < A_k B_k^{-(n+1)},
\end{equation}
where \( u := n + 1, \, v := 2\ell + 1, \, w := n + 1 \). 
Using Python, we verified that
\[
\tau_k \in [1.59, 1.99], \quad \mu_k \in [0.700657, 1.9927] \quad \text{for all odd } k \in [4, 500].
\]
If $-v>0$, then by Equation \eqref{eqn21} it follows that $(n+1)\tau_k<A_kB_k^{-(n+1)}+\mu_k$. However, we have verified that this inequality cannot hold when $n>k^3$ for all odd $k\in[4,500]$. Since we are already operating under the assumption $n>k^3$, it follows that $v>0.$ 
To proceed, we define $M:=3\times 10^{47}$ as an upper bound on $n+1$ in accordance with inequality \eqref{eqn8}. We then apply Lemma \ref{lem1} to Inequality \eqref{eqn21} for each odd value of $k\in[4,500].$ Using Python for computational assistance, we find that within the range
 $M\in [2.8\times 10^{43}, 3\times10^{47}],$ the following estimates hold:
 \begin{align*}
 &1.68\times10^{44}<q_{87} \text{ and } 1.8\times10^{48}<q_{89},\\
 &0.080624\leq \epsilon_k\leq 0.385520,\\
 &1568\leq R_k\leq 130068833,
\end{align*}
where $q_{m_k}$ denotes the denominator of the $m_k^{th}$ convergent in the continued fraction expansion of $\tau_k,$ and $R_k$ represents the upper bound on $n$, obtained for each odd $k$ in $[4, 500].$ Then we recursively apply Lemma \ref{lem1} on the values obtained above and find that the above parameters are in the following ranges:
 \begin{align*}
 &10\leq m_k \leq 21, \hspace{0.2cm} 15009\leq q_{m_k}\leq 135304487,\hspace{0.2cm}
 0.021090\leq \epsilon_k\leq 0.0944463, \hspace{0.2cm} \text{and }\\
 &18\leq R_k\leq 5503681.
\end{align*}
We summarize this result in the following lemma:
\begin{lemma}\label{lem4.2}
Let \( k \in [4, 500] \) be odd and $n_k\geq0$ such that $P_{-n}^{(k)}=0$. Then 
$
n_k\leq R_k\in[18,5503681].
$
\end{lemma}

\subsection*{Derivation of the Zero-Multiplicity Formula}

We define the $k$-generalized Pell sequence $(P_n^{(k)})_{n \in \mathbb{Z}}$ extended to negative indices using the recurrence: $(G_n^{(k)})_{
 n\geq0}$ which has initial values
 $G_i^{(k)} 
 =0$ for $i = 0,1,\ldots,k -2$ and $$G_n^{(k)}=G_{n-k}^{(k)}-
2 G_{n-(k-1)}^{(k)}-\cdots-G_{n-1}^{(k)}$$ holds for all $n \geq k.$
 So, \begin{equation*}
 \begin{split}
G_{n+1}^{(k)} &=G_{n+1-k}^{(k)} - 2G_{n+2-k}^{(k)} - \cdots - G_n^{(k)}\\&= G_{n+1-k}^{(k)} - 2G_{n+2-k}^{(k)} - \cdots -G_{n-1}^{(k)}- \left(G_{n-k}^{(k)} - 2G_{n-(k-1)}^{(k)} - \cdots - G_{n-1}^{(k)} \right)\\&=3G_{n-(k-1)}^{(k)} -  G_{n-k}^{(k)}-G_{n+2-k}^{(k)}, \quad\text{ for all } n\geq k.
\end{split}
\end{equation*}
Therefore
\begin{equation} \label{eq:negative}
G_n^{(k)} = 3G_{n-k}^{(k)} - G_{n - k - 1}^{(k)} - G_{n+1-k}^{(k)}
\end{equation}
is valid for $n \geq k+1$. 
Using the initial conditions of $(G_n^{(k)})_{n\geq 0},$ we observe that $$ G_{n-k}^{(k)}=G_{n - k - 1}^{(k)}=G_{n+1-k}^{(k)}=0 $$
 for all $0 \leq n-k-1 < n-k<n-(k-1)\leq k-2.$ By Equation \eqref{eq:negative}, we get $G_n^{(k)}=0$
 for all $n\in  [k+1,(k-2)+(k-1)].
$
From the above we deduce that  $$G_n^{(k)}=0
 \text{ for all } k+1\leq n-k-1 < n-k<n-(k-1)\leq (k-2)+(k-1) .$$
 Thus, we obtain $G_n^{(k)}=0$
for all $n\in  [2(k+1),2(k-1)+(k-2)].
$
Using computational data and recurrence behavior, one can see that 
the zero set $\Zstroke_k$ appears as a union of consecutive intervals
\[
\Zstroke_k =  \bigcup_{j=1}^{r} I_j, 
\]
where each $I_j$ is a contiguous interval (or ``block") of zeros with $$
r=
\begin{cases}
\vspace{0.2cm}\dfrac{k}{2}, & \text{if } k \text{ is even} \\
\dfrac{(k - 1)}{2}, & \text{if } k \text{ is odd}
\end{cases}
\quad \text{for all } k \geq 4.
$$
The structure of each zero interval $I_j$ is given by
\[
I_j = [(j-1)(k+1), (j-1)(k-1)+(k-2)]
\]
and the total number of zeros in $\Zstroke_k$ is
\begin{equation*}
\chi_k = \sum_{j=1}^{r} j.
\end{equation*}
\vspace{0.3cm}Here, $r$ depends on the parity of $k$:\\
{\bf Case 1: $k$ is Even.}
For $k=2r$, the block lengths are the first 
$r$ odd positive integers in descending order:
$2r-1,2r-3,\ldots,1$
and their sum is
$$\chi_k=(2r-1)+(2r-3)+\cdots+1=r^2.$$
Since 
$k=2r$, so $$\chi_k=\biggr\lfloor\dfrac{k^2}{4}\biggr\rfloor.$$
{\bf Case 2: $k$ is Odd.} Here we put
$k=2r+1$ and the block lengths are the first 
$r$ even positive integers in descending order:
$2r,2(r-1),\ldots,2$, and their sum is
$$\chi_k=2(1+2+\cdots+r)=2\cdot\dfrac{r(r+1)}{2}
=r(r+1).$$
Now compute
$$\biggr\lfloor\dfrac{k^2}{4}\biggr\rfloor=\biggr\lfloor\dfrac{(2r+1)^2}{4}\biggr\rfloor=\biggr\lfloor\dfrac{4r^2+4r+1}{4}\biggr\rfloor=r(r+1),
$$
so again 
$$\chi_k=\biggr\lfloor\dfrac{k^2}{4}\biggr\rfloor.$$
This formula matches all the observed data up to $k = 500$, and the intervals can be recursively constructed using the recurrence.
\subsection*{Conclusion}
We have established that
$
\chi_k =\lfloor k^2/4\rfloor
$
for all $k \in [4,500].$
This completes the proof of Corollary \ref{coro 1.2}.\\

{\bf Data Availability Statements:} Data sharing is not applicable to this article as no datasets were generated or analyzed during the current study.

{\bf Funding:} The authors declare that no funds or grants were received during the preparation of this manuscript.

{\bf Declarations:}

{\bf Conflict of interest:} On behalf of all authors, the corresponding author states that there is no Conflict of interest.
\qed

\end{document}